\documentclass[12pt]{amsart}

 \newtheorem{thm}{Theorem}[section]
 
 \newtheorem{lem}[thm]{Lemma}
 \newtheorem{prop}[thm]{Proposition}
 \theoremstyle{definition}
 
 \theoremstyle{remark}
 
 \theoremstyle{definition}

 \newcommand{\Real}{\mathbb{R}}
 
 \newcommand{\CC}{\mathbb{C}}

 \newcommand{\PP}{\mathbb{P}}

\usepackage{xypic} \usepackage{amsthm}

\begin{document}

\title{Monomials as sums of powers: the Real binary case}

\author[M. Boij]{Mats Boij}
\address[M.Boij]{Department of Mathematics, KTH, Stockholm, Sweden}
\email{boij@kth.se }

\author[E. Carlini]{Enrico Carlini}
\address[E. Carlini]{Dipartimento di Matematica, Politecnico di Torino, Turin, Italy}
\email{enrico.carlini@polito.it}

\author[A.V. Geramita]{Anthony V. Geramita}
\address[A.V. Geramita]{Department of Mathematics and Statistics, Queen's University, Kingston, Ontario, Canada, K7L 3N6 and Dipartimento di Matematica, Universit\`{a} di Genova, Genoa, Italy}
\email{Anthony.Geramita@gmail.com \\ geramita@dima.unige.it  }


\begin{abstract}
We generalize an example, due to Sylvester, and prove that any
monomial of degree $d$ in $\Real[x_0, x_1]$, which is not a power of a
variable, cannot be written as a linear combination of fewer than $d$ powers of linear forms.
\end{abstract}

\maketitle

\section{Introduction}

It is well-known, and easy to prove, that if $k$ is a field of characteristic zero and $R = k[x_0, \ldots , x_n] = \bigoplus_{i=0}^\infty R_i$ is the standard graded polynomial algebra, then the $k$-vector space $R_d$ (for any $d$) has a basis consisting of polynomials $\{ L_1^d, \ldots , L_s^d \}$ where $s = {d+n\choose n} = \dim_kR_d$ and the $L_i$ are pairwise linearly independent forms in $R_1$.  It follows that every form in $R_d$ is a $k$-linear combination of at most $s$ $d^{th}$ powers of linear forms and, if $k$ is algebraically closed, simply a sum of at most $s$ $d^{th}$ powers of linear forms.  We will call such a way of writing $F \in R_d$ a {\it Waring expansion of $F$} because of the echo of Waring's problem from number theory.  We will further refer to such an expression as a {\it minimal Waring expansion for $F$} if the number of summands in such an expression for $F$ is minimal among all such representations.

If $n>0$ and $d=2$ it is a classical fact that although $s = {n+2\choose 2}$ every quadratic form has a Waring expansion involving $\leq n+1 < s$ summands and that, in general, i.e. for $[F]$ belonging to a non-empty Zariski open subset of $\PP(R_2)$ a minimal Waring expansion for $F$ has exactly $n+1$ summands.

These observations have led to a series of problems, usually called {\bf Waring Problems}, which ask for information on minimal Waring expansions for forms of degree $d$ in $R$.

The long outstanding problem of finding the number of summands in a minimal Waring expansion of the generic form of degree $d$ was solved, after being open for almost 100 years, by J. Alexander and A. Hirschowitz (see \cite{AH95}), when $k$ is an algebraically closed field.

Of course, solving this problem for the generic form of degree $d$ does not always give information about any specific form of degree $d$ and the problem of finding the length of the minimal Waring expansion for specific forms has also been a continuing source of interesting speculations and lovely results.  E.g. it was Sylvester (\cite{Harris}) who first observed that although for $R = \CC[x_0,x_1]$, the generic form of degree $d$ has a Waring Expansion with $s = \lceil{{d+1}\over 2}\rceil$ summands, the monomial $x_0x_1^{d-1}$ has $d$ summands in its minimal Waring expansion (the maximum possible).

The Waring problem for specific forms has been considered in depth by B. Reznick in his monograph (see \cite{Reznick}) and by Comas and Seiguer who, to our knowledge, were the first to resolve the problem completely and algorithmically in $\CC[x_0, x_1]$ in their unpublished work (\cite{ComasSeiguer}).

It is interesting to note that although the Waring problem is a very interesting and stimulating problem in purely algebraic terms, it has a surprising number of intimate connections with problems in areas as seemingly disparate as algebraic geometry and communication theory (see for example \cite{RS00},\cite{CaCh} and \cite{ComMour})

Indeed, if $k = \Real$, the field of real numbers, the connection with real world problems is very direct.  This has prompted a re-examination of the Waring problem for $R = \Real[x_0,x_1]$, and a recent very suggestive paper of Comon and Ottaviani (see \cite{ComonOttovaiani}) considered this very problem for degrees $d \leq 5$.

Our main result in this paper follows the line of Sylvester's examples and concerns the minimal Waring expansion for monomials in $\Real[x_0,x_1]$.  We first give a new proof of the fact that the minimal Waring expansion of the monomial $x_0^ax_1^b$ in $\CC[x_0,x_1]$ with $0 <a \leq b$ has $b+1$ summands.  In sharp contrast to this we show that in $\Real[x_0,x_1]$ every monomial of degree $d$ (except $x_0^d$ and $x_1^d$) has $d$ summands in its minimal Waring expansion.

\section{Basic results}

Let $S=k[x_0,x_1]$ and $T=k[y_0,y_1]$. We make
$S$ into a $T$-module using differentiation, i.e. we think of $y_0 = \partial/ \partial x_0$ and $y_1 = \partial/\partial x_1$. We refer to a polynomial in $T$ as $\partial$ instead of using capital letters.  In particular, for
any form $F$ in $S_d$ we define the ideal $F^\perp\subseteq T$ as
follows:
\[F^\perp=\left\{\partial\in T : \partial F=0\right\}.\]

The following {\it Apolarity Lemma} is due to Iliev and Ranestad  \cite{IlievRanestad}.

\begin{lem}\label{apolarityLEMMA}
A homogeneous form $F\in S$ can be written as
\[F(x_0,x_1)=\sum_{i=1}^r \alpha_i(L_i)^d , \ L_i \hbox{ pairwise linearly independent, $\alpha_i \in k$ }\]
i.e. has a Waring expansion with $r$ summands, if
and only if the ideal $F^\perp$ contains the product of $r$
distinct linear forms.
\end{lem}

\section{Binary monomials: the complex case}

The complex case is straightforward for monomials.

\begin{prop}\label{complexsum}
Let $M=x_0^ax_1^b $ be a monomial in $\CC[x_0,x_1]$. If $0<a\leq b$, then $M$ has a minimal Waring expansion with $b+1$ summands, i.e. is a sum of $b+1$ powers of linear forms and no fewer.
\end{prop}
\begin{proof}
Let $I=M^\perp=(y_0^{a+1},y_1^{b+1})$ and notice that the linear
system defined by $I_{b+1}$ is base point free on
$\PP^1=\mathbb{P}S_1$. Applying Bertini's Theorem, we get that the
generic element of $I_{b+1}$ defines a set of $b+1$ distinct
points and hence it is the product of $b+1$ distinct linear forms.
Thus the apolarity lemma yields that $M$ is the sum of $b+1$
powers of linear forms. If $r<b+1$, then $r$ powers do not suffice
as no element in $I_r=(y_0^{a+1})_r$ is a product of $r$ distinct
linear forms.
\end{proof}

\section{Binary monomials: the real case}

We can also ask for a real Waring expansion of a
monomial $M$. More precisely, we want to write

\[M(x_0,x_1)=\sum_{i=1}^r \alpha_i(L_i)^d, \ \ \alpha_i \in \{1,-1 \} \]
where the linear forms $L_i$ are in $\Real[x_0,x_1]$. In
order to do this, we have to increase the number of summands in Proposition \ref{complexsum}.

The following elementary facts will be extremely useful.

\begin{lem}\label{zerocoeff}
Consider the degree $d$ polynomial
\[ F(x)=c_dx^d+\ldots c_1x+c_0\in\Real[x] .\]
If $c_i=c_{i-1}=0$ for some $1 \leq i \leq d$, then $F(x)$ does not have $d$ real roots.
\end{lem}
\begin{proof}  The proof is obvious if $i = 1$ or $i = d$, so we may as well assume that $1 < i < d$.

Consider all the pairs $(c_r,c_s)$ of non-zero coefficients such
that $r>s$ and $c_j=0$ if $r>j>s$. Let $\alpha$ be number of pairs
such that $r-s$ is odd and $\beta$ the number of pairs such that
$r-s$ is even. Notice that, by hypothesis, $\alpha+2\beta<d-1$

Now we apply Descartes' rule of signs. For a pair $(c_r,c_s)$ such
that $r-s$ is odd we get a real root of $F(x)$. For a pair
$(c_r,c_s)$ such that $r-s$ is even we get either two real roots
of $F(x)$ or none.

In conclusion, the number of real roots of $F(x)$ is at most
$\alpha+2\beta$ and we are done.
\end{proof}

\begin{lem}\label{allrealroots}
For each $i<d$ there exists a degree $d$ polynomial
$F(x)=c_dx^d+\ldots c_1x+c_0\in\mathbb{R}[x]$ having $d$ real
roots and such that $c_i=0$.
\end{lem}
\begin{proof}
Choose $a_1,\ldots,a_{d}\in\mathbb{R}$ and consider the polynomial
$F(x)=(x-a_1)\cdot\ldots\cdot (x-a_d)$. This polynomial can also
be written as
\[F(x)=\sum_{i=0}^d E_i(a_1,\ldots,a_d)x^i,\]
where $E_i$ is the degree $i$ elementary symmetric function in its arguments.
The vanishing of the $i$-th coefficient of $F(x)$ can be written
as
\[E_i(a_1,\ldots,a_{d-1})+a_dE_{i-1}(a_1,\ldots,a_{d-1})=0.\]
Hence, if we choose the $a_1,\ldots,a_{d-1}>0$ and distinct there
exists a unique, negative value of $a_d$ such that the coefficient
of $x^i$ in $F(x)$ is zero. As the roots of $F(x)$ are
$a_1,\ldots,a_d$ the polynomial has $d$ real, distinct roots.

\end{proof}

Using the previous results we immediately get a lower bound on the
number of summands in the minimal Waring expansion of a monomial in $\Real[x_0,x_1]$.
\begin{lem}
Let $M=x_0^ax_1^b $ be a monomial in $\Real[x_0,x_1]$. If $0<a\leq b$, then $M$ does not have a Waring expansion with $r\leq a+b-1$ real summands.
\end{lem}
\begin{proof}
Let $I=M^\perp=(y_0^{a+1},y_1^{b+1})$. The general degree $r$
element in $I$ has the form $F(y_0,y_1)=$
\[c_ry_0^r+c_{r-1}y_0^{r-1}y_1+\ldots+c_{a+1}y_0^{a+1}y_1^{r-a-1}+c_{r-b-1}y_0^{r-b-1}y_1^{b+1}+\ldots+c_0y_1^r.\]

If $a+1\geq r-b+2$, then by Lemma \ref{zerocoeff} $F(y_0,y_1)$ is
not the product of $r$ real linear forms. The conclusion follows
by the apolarity lemma.
\end{proof}

\begin{prop}
Let $M=x_0^ax_1^b$ be a monomial in $\Real[x_0,x_1]$. If $0<a\leq b$, then $M$ has a minimal Waring expansion with $a+b$ summands which are powers of real linear forms.
\end{prop}
\begin{proof}
We have that  $M^\perp= I = (y_0^{a+1},y_1^{b+1})$.  Notice that $I_{a+b}$ is
the subspace of $T_{a+b}$ of polynomials which are missing all the
monomials having factor $y_0^a$ or $y_1^b$. Thus, Lemma \ref{allrealroots} and the apolarity lemma yield the result.
\end{proof}

\section{Acknowledgement}  This project started with a research visit
supported by a grant for visiting researchers from the G\"oran Gustafsson Foundation. 

\def\cprime{$'$}

\end{document}